# Bridging knowing and proving in mathematics

# An essay from a didactical perspective

Nicolas Balacheff

CNRS

Laboratoire d'informatique de Grenoble (France)

To Adrien Douady

# 1 An ad hoc epistemology for a didactical gap

## 1.1 The didactical gap

More often than not, the problem of teaching mathematical proof has been addressed almost independently from the teaching of mathematical "content" itself. Some curricula have exposed learners to a significant amount of mathematics without learning about mathematical proof as such (Herbst, 2002, p.288); others teaching mathematical proof as a subject in itself without significantly relating it to concrete practical examples (cf. Usiskin, 2007, p.75). The most common didactical tradition chooses to introduce proof in the context of geometry—usually at the turn of the 8$^{th}$ grade—while completely ignoring it in algebra or arithmetic, where things seem to be reduced to 'mere' computations. This orientation has changed slightly in the past decade with an increasing emphasis on the teaching of proof. However, an implicit distinction between form and content has lead to references to teaching 'mathematical reasoning' (e.g., NCTM standards) or 'deductive reasoning' (e.g., French national programs) instead of mathematical proof as such which would have moved "form" much more to the forefront of the didactical scene.

Nevertheless, it is generally acknowledged that mathematical proof has specific characteristics, among them a formal type of text (the US vocabulary often refers to "formal proof"), a specific organisation and an undisputable robustness once syntactically correct. These characteristics have given mathematics the reputation of having exceptionally stringent practices as compared to other disciplines, practices that are not socially determined but inherent to the nature of mathematics itself.

Hence, the answer to the question: "*Can one learn mathematics without learning what a mathematical proof is and how to build one?*" is "*No*". But now one can observe a *double* didactical gap: (i) mathematical proof creates a rupture between mathematics and other disciplines (even the 'exact sciences') and (ii) a divide in the course of mathematical teaching during the (almost) standard first 12 years of education (into an era before the teaching of proof and one after).

The origin of these gaps lies at the crosspoint of several lines of tension: rigor versus meaning, internal development versus application-oriented development of mathematics, ideal objects defined and manipulated by symbolic representations versus experience-based empirical evidence. I do not analyse these tensions here; I mention them to evoke the complexity of the epistemological and didactical problems which confront us.

One source of the didactical problems is that teaching must take into account the learners' initial understanding and competence: *We can teach only to ones who know…* The learners' existing knowledge often proves resistant, especially because the learners may have proven its efficiency, as in the case of their argumentative skills. In order to overcome this difficulty, teachers organize situations, *mises en scène* and discourses in order to "convince" or "persuade" learners (in the vocabulary of Harel & Sowder, 1998). Argumentation seems the best means to this end. It works both as a tool for teaching and as a tool for doing mathematics for a long while. But then learners suddenly face an unexpected revelation[1]: *In mathematics you don't argue, you prove…*

Looking to bridge this transition, mathematics educators have searched for ideas in psychology. In the middle of the 20th century, the success of Piaget's 'stage theory' of development suggested that proof could be taught only after the required level of development had been reached[2]. As a result, mathematical proof was introduced suddenly in curricula (if at all) in the 9th grade – generally, the year that students have their 13th birthday. However, this strategy has not worked so well, suggesting to some that Piaget may have been wrong.

---

[1] Argumentation means here "verbal, social and rational activity aimed at convincing a reasonable critic of the acceptability of a standpoint by putting forward a constellation of one or more propositions to justify this standpoint" (van Eemeren *et al.*, 2002, p.xii). "In argumentative discussion there is, by definition, an explicit or implicit appeal to reasonableness, but in practice the argumentation can, in all kinds of respects, be lacking of reasonableness. Certain moves can be made in the discussion that are not really helpful to resolving the difference of opinion concerned. Before a well-considered judgment can be given as to the quality of an argumentative discussion, a careful analysis as to be carried out that reveals those aspects of the discourse that are pertinent to making such a judgment concerning it reasonableness." (ibid., p.4)

[2] See e.g. Piaget J. (1969) p. 239: "L'enfant n'est guère capable, avant 10-11 ans, de raisonnement formel, c'est-à-dire de déduction portant sur des données simplement assumées, et non pas sur de vérités observées". More precisely, For more, c.f. Piaget J. (1967) Le jugement et le raisonnement chez l'enfant. Delachaux et Niestlé.

Some mathematics educators then turned to psychologies of discourse and learning, feeling that the followers of Piaget had not paid enough attention to language and social interaction. Some suggested the ideas of Vygotsky and the socio-constructivists could have provided a solution (e.g. Forman *et al*. 1996). However, this line of thought did not appear to be the panacea either. Then Lakatos' work seemed to suggest that a solution might be found in the epistemology of mathematics itself (e.g. Reichel 2002); however, such attempts also failed amid scepticism from mathematicians and researchers.

The responsibility for all these failures does not belong to the theories which supposedly underlie the educational designs, but to naive or simplifying readers who have assumed that concepts and models from psychology can be freely transferred to education. In particular, they rarely take into account the nature of mathematics as content (while often emphasizing the nature of the perceived practice of mathematicians).

My objective here is then to question the constraints mathematics imposes on teaching and learning, postulating that, as for any other domain, learning and understanding mathematics cannot be separated from understanding its intrinsic means for validation: *mathematical proof*. First, I address the epistemology of proof, on which we could base our efforts to manage or bridge the didactical gap discussed above.

## *1.2 The need to revisit the epistemology of proof*

Although apparently a bit simplistic, it may be good to start from the recognition that mathematical ideas are not a matter of feeling, opinion or belief. They are of the order of 'knowing' in the Popperian sense[3], by virtue of their very specific relation to proof (and proving). They provide tools to address concrete, materialistic or social problems, but they are not about the "real" world. To some extent, mathematical ideas are about mathematical ideas; they exist in a closed 'world' difficult to accept but difficult to escape. For this reason, mathematical ideas do not exist as plain facts but as statements which are accepted only once they have been proved explicitly; before that, they cannot be[4] instrumental either within mathematics or for any application.

However, despite this emphasis on the key role of proof in mathematics, it must be remembered that at stake is not *truth* but the *validity of a statement within a well-defined theoretical context* (cf. Habermas, 1999). For example, Euclidean geometry is no truer than Riemannian geometry. This shift from the vocabulary of truth to the vocabulary of validity, which suggests a shift from *proof* to *validation,* is more important than we may have realized. Validation refers to constructing reasons to accept a specific statement, within an accepted framework shaped by accepted rules and other previously accepted statements. From this perspective, mathematical validation searches for an *absolute proof in an explicit context*; it can thus claim certainty as a foundational principle.

This view of validity and proof is antiauthoritarian (Hanna & Janke, 1996, p.891), insofar as it assumes a common agreement about a collective and well-understood effort. It thus fits the classical conception of what a scientific proof should be, since such a proof must clearly not depend on specific individual or social interests. In this lies the democratic aspect of mathematical proof, as noted by Hanna. Proving is an example of an intellectual enterprise that allows a minority to overcome the opinion of an established majority, according to shared rules. This is related to an ancient meaning of the word "demonstration" in English (e.g., Herbst, 2002, p.287).

---

[3] Popper (1959) proposed falsification as the the empirical criterion of demarcation of knowledge, scientific theories or models.

[4] Or should not be...

So the concept of *proof* is not a stand-alone concept; it goes with the concepts of "validity of a *statement*" and "*theory*". This has been well explained and illustrated by the Italian school, especially Alessandra Mariotti (1997). However, the word "theory" is the most difficult for learners. No such thing is available to learners *a priori*, and to understand what the word means seems out of reach. Nevertheless, learners have ideas about mathematics and about mathematical facts. They also have experience in arguing about the "truth" of a claim or the "falsity" of a statement they reject; but this is experience in argumentation in contexts that are not framed by a theory in scientific terms. To construct a proof requires an essential shift in the learner's epistemological position: passing from a practical position (ruled by a kind of logic of practice) to a theoretical position (ruled by the intrinsic specificity of a theory).

In addition, we cannot engage in the validation of 'anything' that has not been first expressed in a language. This principle applies across disciplines (Habermas, 1999), but plays a special role in mathematics, where the access to 'mathematical objects' depends in the first place on their semiotic availability (Duval, 1995).

In other words, the teaching and learning of mathematical proof requires mastery of the relationships among knowing, representing and proving mathematically.

# 2 A model to bridge knowing and proving

## 2.1 Short story 1: Fabien and Isabelle misunderstanding

--------------------------------------------------------------------------------

Consider the following problem[5]:

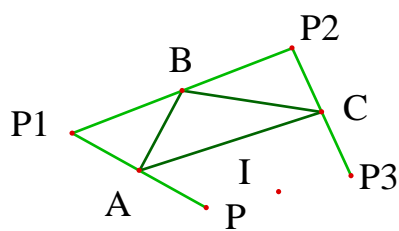

>Construct a triangle ABC. Construct a point P and its symmetrical point P1 about A. Construct the symmetrical point P2 of P about B, construct the symmetrical point P3 of P about C. Move P. What can be said about the figure when P3 and P are coincident? Construct the point I, the midpoint of [PP3]. What can be said about the point I when P is moved? Explain.

Figure 1.

--------------------------------------------------------------------------------

Constructing the diagram (Fig.1) with dynamic geometry software[6], one can easily notice that the point I does not move when one manipulates the point P. This *fact* seems surprising; the crux of the situation is to propose an explanation.

Let us examine the interaction between a teacher and a student, Fabien, about this problem[7]. Fabien has observed the fact but he has no insight about the reason: "*The*

---

[5] From Capponi, 1995, *Cabri-classe*, sheet 4-10.

[6] e.g. Cabri-geometry (here used for the drawing), or Geometer Sketchpad; or Geogebra or one of the several others now available sometimes open access.

[7] A more detailed analysis can be found in Balacheff & Soury-Lavergne (1995), Sutherland & Balacheff (1999).

*point I does not move, but so what...*" However, he noticed and proved that ABCI is a parallelogram. At this stage, from the point of view of geometry (and of the tutor), the reason I stands immobile while P moves should be obvious. The tutor then provides Fabien with several hints but with no results. After a while she desperately insists: *"The others, they do not move. You see what I mean? Then how could you define the point I, finally, without using the points P, P1, P2, P3?"* Throughout the interaction, the tutor is moved by one concern which can be summarized by the question: *"Don't you see what I see?"* But Fabien does not see the 'obvious'; it is only when she tells him the mathematical reasons for the immobility of I that the tutor provokes a genuine "Aha!" effect...

In order to explain the immobility of I, the teacher had get the student to construct a link between a *mechanical world*—that of the interface of the software[8] – and a *theoretical world*— the world of geometry. Only this link can turn the observed *fact* (the immobility of I) into a *phenomenon* (the invariance of I). But the construction of this link is not straightforward; it is a *process of modelling.*

Teacher and student did share representations, words, and arguments so that they could communicate and collaborate; however, this did not guarantee that they shared understanding. Educators have made considerable efforts to develop representations that could make the nature and the properties of mathematical concepts more tangible. But these remain just representations with no visible referent; manipulating them and sharing factual experience does not guarantee shared meaning. Nevertheless, they are the only means of communication, since in mathematics the referent, in a semiotic sense, is itself a representation (i.e., a tangible entity produced on purpose).

---

[8] Another student's search for an explanation illustrates well what is meant here by mechanical world: "*So... I have said... But is not very clear... That when, for example, we put P to the left, then P3 compensates to the right. If it goes up, then the other goes down...*" (Sébatien, [prot. 78-84]).

In the next section, I will explore this issue of representation and its relation with the learners' building of meaning, and then take up the challenge of defining "knowing" in a way that may not solve the old epistemological problem but will provide some grounds to build a link between knowing and proving.

## *2.2 Trust, doubt and representations*

The fascination for proof without words[9], which would give access to the very meaning of the validity of a mathematical statement without the burden of sophisticated and complicated discourses, is a symptom of the expectations mathematics educators have attached to the use of non-verbal representations in mathematics teaching. The development of multimedia software, advanced graphical interfaces and access to 'direct manipulation' of the represented 'mathematical objects' has even strengthened these expectations. The above story of the Fabien and his tutor misunderstandings is initial evidence that things might be slightly more difficult. I will explore this difficulty now, starting with an example coming from professional mathematics.

In 1979, Benoit Mandelbrot noticed in a picture produced by a computer and a printer that the Mandelbrot set[10]—as it is now known, following a suggestion of Adrien Douady—was not connected. "A striking fact, which I think is new" Mandelbrot[11] remarked. John Hubbard, a former PhD student of Adrien Douady's who became his well known collaborator, reported that:

---

[9] See Claudi Alsina and Roger B. Nelsen (2006), *Math Made Visual: Creating Images for Understanding Mathematics*, published by MAA, and a good example in Roger B. Nelsen (1993), *Proofs without words: exercises in visual thinking*, published by MAA. See Hanna (2000, esp. pp.15-18) for an analysis.

[10] Considering the sequence of complex numbers $z_{n+1} = z_n^2 + c$, the Mandelbrot set (or set M) is obtained by fixing $z_0=0$ and varying the complex parameter c.

[11] Quotation from p.250 of Mendelbrot (1980) Fractal aspects of the iteration of $z \rightarrow \lambda z(1-z)$ for complex $\lambda$ and z. Annals of the New York Academy of Sciences. 357 (1) 249 - 259

> Mandelbrot had sent [them] a copy of his paper, in which he announced the appearance of islands off the mainland of the Mandelbrot set M. Incidentally, these islands were in fact not there in the published paper: apparently the printer had taken them for dirt on the originals and erased them. (At that time, a printer was a human being, not a machine). Mandelbrot had penciled them in, more or less randomly, in the copy [they] had. (Hubbard 2000 pp.3-4)

This anecdote reflects two things: first, the efficiency and strength of the computer-based picture in supporting a conjecture; second, the fragility of this same picture, which depends on both the algorithmic and technical conditions of its production. Then, Hubbard reported:

> One afternoon, Douady and I had been looking at this picture, and wondering what happened to the image of the critical point by a high iterate of the polynomial $z^2 + c$ as c takes a walk around an island. This was difficult to imagine, and we had started to suspect that there should be filaments of M connecting the islands to the mainland. (ibid.)

Soon, Adrien Douady realized that this meant that the set M is connected[12], but "the proof of this fact is by no means obvious," he remarked (Douady, 1986, p.162). The proof followed after a long process of writing, initiated by a *Note aux Comptes-rendus* in 1982. After the discovery of the connectedness, images of the set M got transformed, offering a more beautiful picture full of colours which, so to speak, 'displayed' the connectivity of M (Fig. 2).

---

[12] Régine Douady remembers that Adrien had been quickly convinced of the connectivity of M, thanks to the theoretical argument which convinced him in an astonishingly "simple" way. However, to complete the explicit proof took some time (*2008, personal communication*).

---

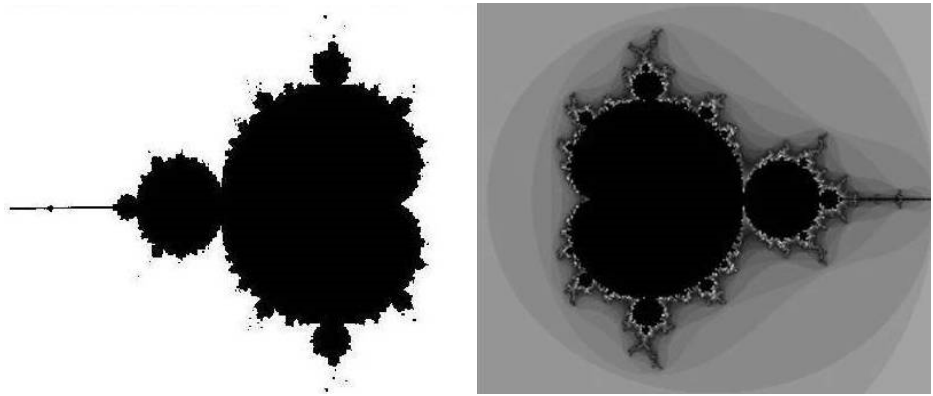

The Mandelbrot set for z→z$^2$+c

before and after the Douady and Hubbard discovery

Figure 2.

---

This case supports the idea of complex relations between representation and mathematical objects—or, more precisely, the role of representations as mediators for the conceptualisation of mathematical objects. It invites more caution in considering evidence in a non-verbal representation. Not to say that non-verbal representations or expressions of an argument are of no value; rather, I emphasize that the frequent claim in education that, *"A picture is worth a thousand words"* has limits and cannot be accepted without further examination.

For example, graphic calculators are widely used by students. They provide students with efficient tools for calculus, blending graphical and symbolic representations. The use of this technology has led to new problem-solving strategies that take advantage of the low cost of exploring of graphical representations. Among them is what Joel Hillel (1993, p.29) called "window shopping," which consists of playing with the various possibilities offered by the display. The diagrams (below) reproduce two appearances of the graph of the same function, $f(x) = x^4-5x^2+x+4$. As one can 'see', these pictures can induce different conjectures about, for example, the numbers of zeros of the polynomial or its behaviour within the interval [-2, +2]

---

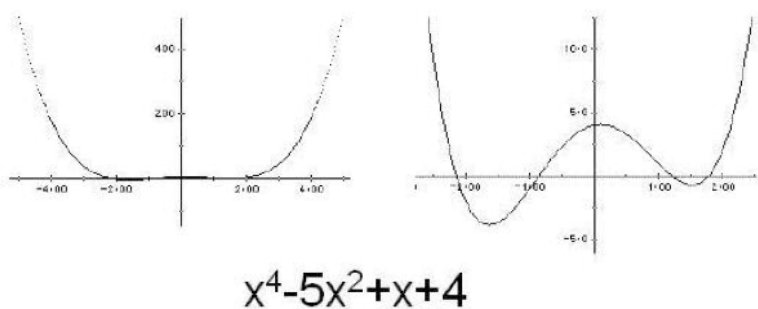

Figure 3.

---

It is now common for teachers to warn students and teach them strategies to ensure reliable, optimal use of their calculators. Still, the problem of knowing how to balance trust and doubt when using these machines and looking for conjectures has no straightforward answer. Part of achieving this balance depends not only on how learners critically organize their explorations but also on the reliability of the embedded software. Consider the case of the function $g(x)=\sin(e^x)$. Most students are prepared to study this function without a priori foreseeing difficulties; that is, until their machine displays something like the following picture:

---------------------------------------------------------------------------------------------------

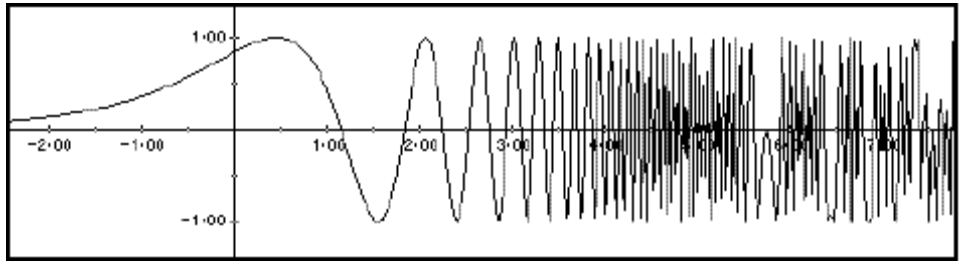

Figure 4.

---------------------------------------------------------------------------------------------------

'Window shopping' will not help to answer the questions this display raises. An algebraic study will just leave students with a question they probably cannot solve with their knowledge of mathematics and computer science. This picture results from the interference between the computation of the coordinates of each point to be displayed and the choice of which pixel to turn black on the screen. In the end, it is the product of a kind of stroboscopic effect, as suggested by Adrien Douady[13]. Producing a 'correct' figure would be a matter of first mathematically notating both the capabilities and the limitations of the drawing instrument and then using sophisticated computational strategies to decide on the intervals at which to plot an 'informative' graph.

The problem of how students can decide to trust or doubt mathematical representations goes beyond graphical representations to include any representation. A last example, taken from Luc Trouche work (2003) on computer algebra systems demonstrates this. Consider the equation $[Ln(e^x-1)=x]$: One can use a pocket graphical calculator to solve it algebraically or to graph it; the two pictures below (from Trouche, 2003, p.27) show the respective results.

---

[13] Personal communication

-------------------------------------------------------------------------------------------------

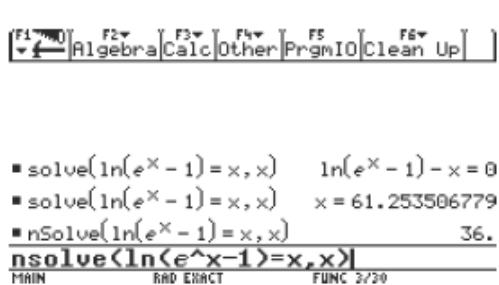 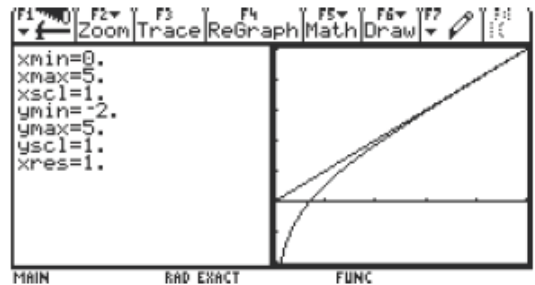

Figure a. *Traces de l'étude algébrique et numérique*   Figure b. *Traces de l'étude graphique*

Figure 5

-------------------------------------------------------------------------------------------------

The results speak for themselves. The optimal treatment leading to a solution – in this case, that this equation has no solution – consists of a formal transformation of the algebraic expression, producing [$e^x-1= e^x$].

The difficulty students may have relates not to their lack of mathematical knowledge but to a general human inclination not to question their knowledge and their environment unless there is a tangible contradiction between what is expected after a given action and what is obtained, as my final example will demonstrate.

In this case, upper secondary students were asked to tell what is the limit at $+\infty$ of the function [$f(x)=\ln(x)+10\sin(x)$]. Without a graphic calculator, only five percent of the students answered wrongly; with a graphic calculator, which displayed the window reproduced below, this number grew to 25 percent (Guin & Trouch, 2001, p.65).

---

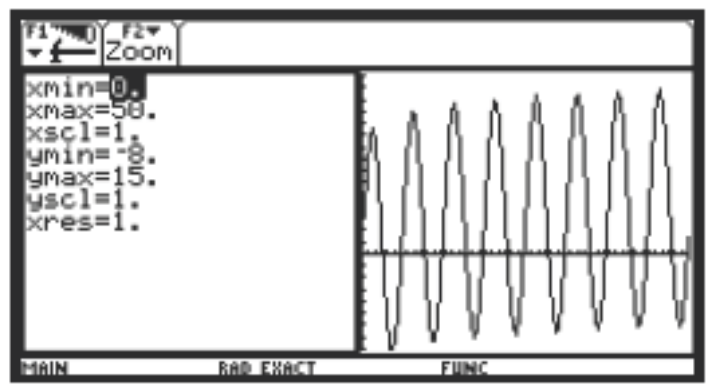

Figure 6.

---

Given such cases of error, teachers and mathematics educators might have to consider whether graphic calculators contribute positively to mathematics learning or whether students have difficulty shifting from one semiotic context to another. (Other examples of common errors include: the value of π is exactly 3.14, or a convergent series reaches its limit, or the Fibonacci series $U_0=1$, $U_1=(1+\sqrt{5})/2$, $U_n=U_{n-1}+U_{n-2}$ is divergent). Most such errors, or "misconceptions" to use the 1980s term, are probably symptomatic of the students' knowledge, which can be legitimate in certain contexts although possibly wrong mathematically. To analyse this issue further, we must have a conceptualization of the students' knowledge which (i) allows us to make sense of it from a mathematical perspective; (ii) is relevant from a cognitive perspective; and (iii) opens the possibility of didactical solutions.

## 2.3 A phenomenological definition of knowing

Studying students' productions that were mathematically incorrect, the mathematics educators of the 1980s usually chose to use the word "misconception". As noted by Jere Confrey (1990), such student errors should be first considered as indications of what they know. Comfrey used the generic word "conception" to refer to the rationale of students' answers to a given problem or question. I postulate that such conceptions result from the learner's interactions with the environment, and that learning is both a process and an outcome of the learner's adaptation to this environment. By "environment", I refer to a physical setting, a social context or even a symbolic system (especially now that the latter can be depicted by a technology which dynamically materialises it).

However, only some characteristics of the environment are relevant from the point of view of learning. Educators do not deal with the learner in all his or her social, emotional, physiological and psychological complexity, but from a knowledge perspective: as *the epistemic subject*. The same principle applies to the environment, which we restrict to *the milieu* defined as *the subject's antagonist system* in the learning process (Brousseau, 1997, p.57); that is, we only consider those features of the environment that are relevant from the epistemic perspective. This means that our characterizations of the (epistemic) subject and of the milieu are interdependent systemically (and dynamically, since both will evolve during the learning process).

Pragmatically, the only accessible evidences of a conception are behaviours and their outcomes. Our problem is to interpret these in terms of indicators of strategies the adapted nature of which must be demonstrated in a model or representation attributed to the student (Brousseau, 1997, p.215)[14]. The formalisation of a conception I propose below aims at providing such a model. Recognizing this interdependence, expressed by Noss and Hoyles[15] (1996, p.122) as *situated abstraction*, accepts that people could demonstrate different and possibly contradictory conceptions depending on circumstances, although knowledgeable observers may ascribe them to the same source concept.

Thus, a conception is attached neither to the subject nor to the milieu, but exists as *a property of the interaction between the subject and the milieu*—its antagonist system (Brousseau, 1997, p.57). The objective of this interaction is to maintain the viability of the subject/milieu system (or [S↔M] system) by returning it to a safe equilibrium after some perturbation (i.e., the tangible materialization of a problem). This implies that the subject recognizes the perturbation (e.g., a contradiction or uncertainty) and that the milieu has features which make the perturbation tangible (since otherwise, the milieu may "absorb" or "tolerate" errors or dysfunctions).

---

[14] For the convenience of the English-speaking reader, I take all the references to Brousseau's contributions to mathematics education from Kluwer, 1997 but Brousseau's work was primarily published between 1970 and 1990.
[15] This proposition should be understood in the light of the development of the 'situated learning paradigm' of Jeane Lave and Etienne Wenger, whose work was published in the early 1990s.

---

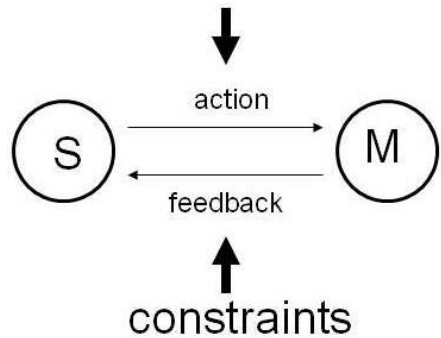

*A* conception *is the state of dynamical equilibrium of an action/feedback loop between a subject and a milieu under proscriptive constraints of viability.*[16]

Figure 7

---

From this definition of *conception*, I can derive a definition of *knowing* as the characterization of a dynamic set of conceptions. This definition has the advantage of being in line with our usual use of the word "knowing" while providing grounds to understand the possible contradictions evidenced by learners' behaviours and their variable mathematical development. A conception is a situated knowing; in other words, it is the instantiation of a knowing in a specific situation detailed by the properties of the milieu and the constraints on the relations (action/feedback) between this milieu and the subject.

---

[16] These constraints do not address how the equilibrium is recovered but the criteria of this equilibrium. Following Stewart (1994, pp. 25-26), I argue that these constraints are *proscriptive* – they express necessary conditions to ensure the system's viability – and not prescriptive, since they do not tell in detail how equilibrium must be reconstructed.

This definition of conception provides a starting point but still has to be refined in order to make it relevant to our research. To do so, I will now introduce the model cK¢[17], in order to provide an effective tool to concretely represent and analyze the corpus of data obtained from the observation of students' activities. This model aims to establish a necessary bridge between knowing and proving by providing a more balanced role to control structures with respect to the role usually assigned to actions and representations.

## 2.4  A model to bridge knowing and proving:  cK¢

That validation plays a key role in the emergence of 'knowing' has been established at least since Popper proposed the criterion of falsification and Piaget introduced the process of cognitive disequilibrium. This principle is also inherent in a "conception" as we define it, adding the explicit condition that a conception is not self-contradictory.

"Proving" is the most visible part of the intellectual activity related to validation. However, as the Italian school has clearly demonstrated (Boero *et al*. 1996), proving cannot be separated from the on-going controlling activity involved in solving a problem or achieving a task. To some extent, "proving" can be seen as an ultimate achievement of controlling and validating. No one can claim to know without a commitment to and a responsibility for the validity of the claimed knowledge.  In return, this knowledge functions as a means to establish the validity of a decision in the course of performing a task and even in the process of building new knowledge—especially in the learning process. In this sense, knowing and proving are tightly related.  Hence, *a conception is validation dependent*: In other words, we can diagnose the existence of a conception because there is an observable domain in which "it works", in which there are means to validate it and to challenge possible falsifications. This is the essence of Vergnaud's (1981, p.220) statement that problems are the sources and criteria of concepts.

---

[17] The letters cK¢ stand for : "conception", "knowing" and "concept"; more about this model is presented and discussed on [http://ckc.imag.fr]

Vergnaud demonstrated that we could characterize students' conceptions with three components: problems, representation systems and invariant operators (1991, p.145)[18]. I take this model as a starting point, with the addition of the related control structure.

Then, I can characterize a conception by a quadruplet (P, R, L, $\sum$) in which:

- P is a set of problems;

> This set corresponds to the class of the disequilibria the considered subject/milieu [S↔M] system can recognize; in mathematical terms: P is the set of problems which can be solved—in pragmatic terms, P is the conception's *sphere of practice*.

- R is a set of operators;
- L is a representation system;

> R and L describe the feedback loop relating the subject and the milieu, namely the actions, feedbacks and outcomes.

- $\sum$ is a control structure;

> The control structure describes the components that support the monitoring of the equilibrium of the [S↔M] system. This structure ensures the conception's coherence; it includes the tools needed to take

---

[18] Vergnaud in fact proposed this definition at the beginning of the 1980's.

> decisions, make choices, and express judgement on the use of an operator or on the state of a problem (i.e., solved or not).

This model aims at accounting for the [S↔M] system and is not restricted to one of its components[19]. The representation system allows the formulation and the manipulation of the operators by the active subject as well as by the reactive milieu. The control structure allows expression and discussion of the subject's means for deciding the adequacy and validity of his or her action as well as the milieu's criteria for selecting a feedback. This symmetry allows us both to take the subject's perspective when evaluating his or her knowing and the milieu's perspective when designing the best conditions to stimulate and support learning. Moreover, it gives us a framework in which to describe, analyze and understand the didactical complexity of learning proof by taking into account the interrelated relevant dimensions: the subject, the milieu and the problem.

In the next section I will give an illustration of this distinctive role of the control structure and the light it sheds on the learners' behaviors we observe and aim at understanding. I will then summarize the proposed framework discussing the relations we must establish between action, formulation and validation in order to understand the didactical complexity of learning and teaching mathematical proof. These three dimensions provide the means we need to build a bridge between knowing and proving.

---

[19] By extension, one can often refer to students' conceptions as acceptable given that one can account precisely for the circumstances, which are the milieu and the constraints within which [S↔M] functioned.

# 3   Proving from a learning perspective

## *3.1  Short story 2: Vincent and Ludovic mismatch*

Vincent and Ludovic are two middle school students who had no specific difficulties with mathematics. They volunteered to participate in an experiment that Bettina Pedemonte (2002) was carrying out to study the cognitive unity between problem solving and proof. The problem was the following:

---

> *Construct a circle with AB as a diameter.  Split AB in two equal parts, AC and CB.  Then construct the two circles of diameter AC and CB… and so on.*

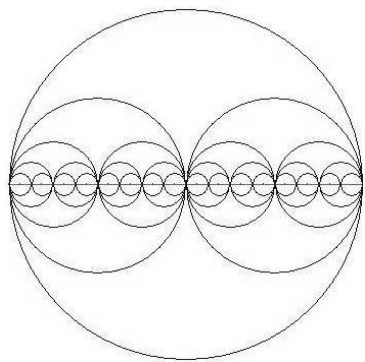

> *How does the perimeter vary at each stage?*
>
> *How does the area vary?*

Figure 8.

---

With no hesitation, the two students expressed – with the formulas they knew well – the perimeter and the area of the first steps in the series of drawings.  Their letters represent quantities and the formulas are another description of the reality the drawing factually displays.  The students conjectured that the perimeter will be constant and that the area decreases to zero. But Vincent noticed that "*the area is always divided by 2...so, at the limit? The limit is a line, the segment from which we started ...*". The discussion then continued:

> 41. *Vincent*: It falls in the segment… the circle are so small.
>
> 42. *Ludovic*: Hmm… but it is always 2πr.
>
> 43. *Vincent*: Yes, but when the area tends to 0 it will be almost equal…
>
> 44. *Ludovic*: No, I don't think so.
>
> 45. *Vincent*: If the area tends to 0, then the perimeter also… I don't know…
>
> 46. *Ludovic*: I will finish writing the proof.

Although Vincent and Ludovic collaborate well and seem to share the mathematics involved, the types of control they have on their problem-solving activity differ. Ludovic is working in the algebraic setting (c.f., Douady, 1985); the control is provided by his ensuring the correctness of the symbolic manipulation and his knowledge of elementary algebra. Vincent is working in a symbolic-arithmetic setting; the control comes from a constant confrontation between what the formula "tells" and what is displayed in the drawings. Both students understood the initial situation in the "same" way, both syntactically manipulated the symbolic representations (i.e., the formulas of the perimeter and of the area), but their controls on what they performed were different, revealing that the conceptions they mobilized were also significantly different. I deduce that the operators they manipulated (algebraic writings, sketching diagrams, etc.), although they coincided from the behavioural perspective, were semantically different. Moreover, from this evidence, an observer could argue that the students were not addressing the same "problem"; Vincent was "baffled" by the gap between what he saw and what he computed, while Ludovic was "blind" to this gap. (Actually, Ludovic's knowledge of calculus would not have been sufficient to provide any relevant explanation).

The symbolic representation plays the role of a semiotic mediator between the two students' different conceptions. It allows communication between the students and is instrumental for each in controlling the problem-solving process and building a proof. We know that two different representations may demonstrate two different understandings; however, here one given representation also supports different understandings and hence different proofs.

## *3.2 The complex nature of proof*

Many theorists have attempted to answer the question of what counts as a proof, from either an epistemological or an educational point of view. However, there is no single, final answer. The Vincent and Ludovic discussion above confirms that sheer formal computation is not enough. As in one of the best previous anecdotes in the history of mathematics[20], Vincent could well say to Ludovic: *I see it, but I don't believe it*. As several authors have emphasised, a proof should be able to fulfil the need for an explanation; however the explanatory nature of a proof may become the object of an even more irreconcilable disagreement than was its rigor. Consider the simple mathematical statement: The sum of two even numbers is itself even. The following figures provide a sample of proofs of this statement. A discussion of these proofs by mathematicians, mathematics teachers and learners provokes very different responses from each.

---

[20] "*Je le vois, mais je ne le crois pas*", wrote Cantor to Dedeking, in 1877, after having proved that for any integer n, there exists a bijection between the points on the unit line segment and all of the points in an n-dimensional space.

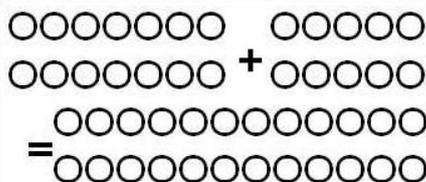

Example adapted from Healey and Hoyles, (2000 p.400)

Figure 9

The arguments in such a discussion involve three types of critical considerations: the search for certainty, the search for understanding and the requirements for a successful communication. The complex nature of proof lies in the fact that any effort to improve a candidate-proof on one of these dimensions may change its value on the other two. There is no clear standard to decide on the correct balance. Restricting the evaluation to the "certainty" side is playing safe, as this side is compulsory for the transformation of mathematical ideas. However, such reductionism is not viable from a learning perspective, especially when students are first introduced to mathematical proof; their control structures are not appropriately evolved. Educators at this point need to give academic status to activities that may not lead to what would be a proof for professional mathematicians but that still make sense as mathematical activities. Hence, my proposal to structure the relations between explanation, proof and mathematical proof as I did to ground my own work (Balacheff, 1988). This structuration distinguished between pragmatic and intellectual proof, and within both it identified categories related first to the nature of the student knowing and his or her available means of representation.

The rationale for this organisation (sketched below in figure 10) is the postulate that the explaining power of a text (or non-textual "discourse") is directly related to the quality and density of its roots in the learner's (or even mathematician's) knowing. What is produced first is an *"explanation" of the validity of a statement from the subject's own perspective*. This text can achieve the *status of proof if it gets enough support from a community* that accepts and values it as such. Finally, it can be claimed as *mathematical proof if it meets the current standards of mathematical practice*. So, the keystone of a *problématiques* of proof in mathematics (and possibly any field) is the nature of the relation between the subject's knowing and what is involved in the 'proof'.

---

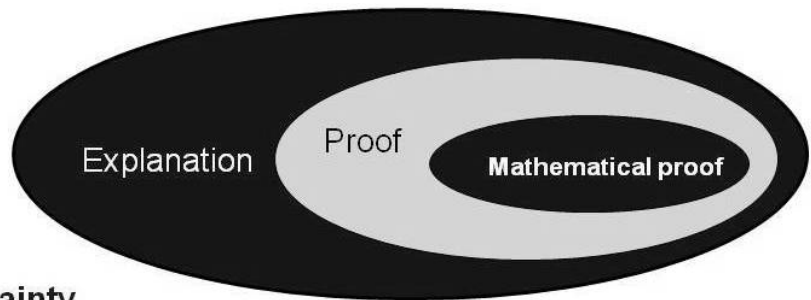

Figure 10

---

This recognition of a proof's roots in knowing may justify a statement as strong as Harel & Sowder's that "one's proof scheme is idiosyncratic and may vary from field to field, and even within mathematics itself," (1998, p.275). However, this view misses the social dimension of proof, which transcends an entirely subjective feeling of understanding (as well as "ascertaining" or "persuading"; Harel & Sowder, ibid., p. 242). From a didactical perspective, the issue is not psychological but epistemological, being directly related to the role a proof plays in building links between a theory that provides its framework and means and a statement that it aims to validate. The transcendence of a proof, proposed by Habermas (1999) as a requirement for a *problématique* of truth and justification, is a dimension too often forgotten in favour of a psychological or sociological analysis of proving. This transcendence is not a dogmatic but a pragmatic position which allows the construction of knowledge as a collective asset which can be shared and be sustainable without depending on its author(s) and circumstance(s) of birth.

The technicalities of mathematical proof are then essential, and can be accepted as the price for a viable construction of mathematics. In this respect, formal rigour is a weapon against the biases that "idiosyncratic proof schemes" may produce.

### *3.3 Knowing and proving in the didactical genesis of proof*

Learning mathematics starts with the first years of schooling, at least from an institutional point of view. As is well documented, learners at this elementary level depend as much on their experience as on the teacher as a reference to distinguish between their opinions, their beliefs and their actual knowledge. The criterion for assessing this difference rests either in the tangible efficiency of the knowledge at stake or in ad hoc validation by the teacher. But the teacher has to rely on knowledge, demonstrating that authority is not the ultimate reference. Hence, efficiency and tangible evidence are the supports for the validity of a statement: It's true because we verify that it works. Mathematical learners are first of all practical persons; to enter mathematics they have to change their intellectual posture and become a theoretician. This shift can easily be seen in the passage from practical geometry (the geometry of drawings and shapes) to theoretical geometry (the deductive or axiomatic geometry), or from symbolic arithmetic (computation of quantities using letters) to algebra. A learner making the transition from the practical to the theoretical has to face the epistemological difficulty of a transition from knowing in action to knowing in discourse: The origin of knowing is in action but the achievement of mathematical proof is in language (see below figure 12).

Again, the tight relationship among action, formulation (semiotic system) and validation (control structure) imposes itself (Brousseau 1997). This trilogy which defines a conception, also shapes didactical situations[21]; there is no validation possible if a claim has not been explicitly expressed and shared; and there is no representation without a semantic which emerges from the activity (i.e., from the interaction of the learner with the mathematical milieu).

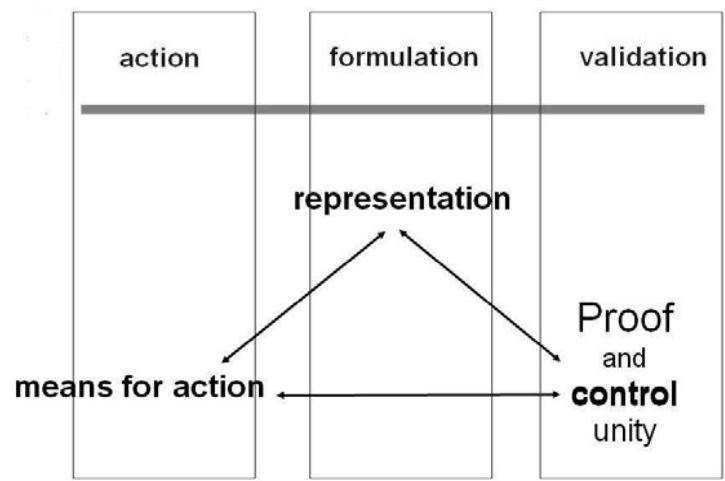

Indeed, this passage from mathematics as a tool whose rationale is 'transparent', to mathematics as a theoretically-grounded means for the production and evaluation of explicit validation has a key stepping stone: language; as a *symbolic technology* (Bishop 1991 p.82), not just a means for social interaction and communication. Language allows learners to understand and appropriate the value of mathematical proof compared with the pragmatic proof they were used to. Now, this language could be of lower levels than the naïve formalism mathematicians use; the level of language will bind the level of the proof learners can produce and/or understand. However, there is room for genuine mathematical activity at all these levels, provided that the learners have moved beyond empiricism and have seen the added value of the theoretical posture (see figure 12 below).

---

[21] figure 11 below sketches the interactions between these three poles

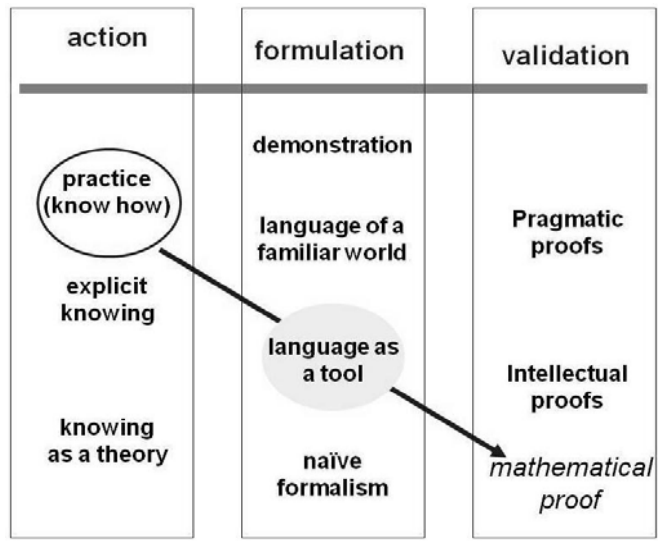

This figure illustrates the approximate mapping between the critical categories in each of the three dimensions (action, formulation and validation). It points the most difficult problem for teacher, that is to provide students with the means to switch from a pragmatic approach of truth to a theoretical approach of validity based on mathematical proof. Realising that language is a tool is a critical milestone on this move.

Figure 12

---

## 4   Still an open problem: the situations…

After a few decades, researchers have now reached a consensus on the variety of meanings that proof may have for learners (if not for teachers). Several classifications and analyses of the complexity of the different aspects of mathematical proof have been extensively reported.  Although they still express significant differences (Balacheff, 2008), researchers have converged on considering mathematical proof as a core issue in the challenge of learning and teaching mathematics; mathematical knowing and proving cannot be separated. In other words an educational *problématique* of proof cannot be separated from that of constructing mathematical knowledge.

This challenge is well understood from an epistemological perspective. However, it is far from clear from a didactical perspective. A lot of effort has gone into proposing problems and mathematical activities which could facilitate the learning of mathematical proof. At the turn of the 20$^{th}$ century, computer science and human-computer interaction research have made so much progress that it is possible to provide learners and teachers with environments able to provide much more mathematically relevant feedback on users' activities. Especially, dynamic geometry environments and computer algebra systems allow learners to experience conjecturing and refuting in a manner never available before, hence giving them access to a dialectic necessary to ground the learning of mathematical proof. However, there is some evidence that learners can remain in a pragmatic intellectual posture, not catching the value of mathematical proof.

Prompting the ultimate move from pragmatic to theoretic knowing requires designing situations so that the pragmatic posture is no longer safe or economical for the learners, while the theoretical posture demonstrates all its advantages. The resultant social and situational challenges are levers which one can use to modify the nature of the learners' commitment to proving. Such design is possible if solving a problem is no longer the main issue and fades away behind the issue of being "sure" of the validity of the solution. We already have some examples which witness the possibility of designing such situations (e.g., Bartolini-Bussi 1996, Boero *et al*. 1996, Arsac and Mantes 1997, etc.). The scientific challenge is now to better understand the didactical characteristics of these situations and to propose a reliable model for their design, for the sake of both researchers and teachers.